\documentclass[12pt]{article}%
\usepackage{graphicx}
\usepackage[intlimits]{amsmath}
\usepackage{latexsym}
\usepackage{amsfonts}
\usepackage{amssymb}%
 \makeatletter
\newfont{\footsc}{cmcsc10 at 8truept}
\newfont{\footbf}{cmbx10 at 8truept}
\newfont{\footrm}{cmr10 at 10truept}
\newtheorem{theorem}{Theorem}

\newtheorem{lemma}[theorem]{Lemma}
\newtheorem{fact}[theorem]{Fact}

\newenvironment{proof}[1][Proof]{\noindent{\textbf {#1}  }}  {\hfill$\Box$\bigskip}

\begin{document}

\title{Tur\'{a}n's theorem inverted}
\author{Vladimir Nikiforov\\{\small Department of Mathematical Sciences, University of Memphis, Memphis TN
38152}\\{\small email: vnikifrv@memphis.edu}}
\maketitle

\begin{abstract}
Let $K_{r}^{+}\left(  s_{1},\ldots,s_{r}\right)  $ be the complete $r$-partite
graph with parts of size $s_{1}\geq2,s_{2},\ldots,s_{r}$ with an edge added to
the first part. Letting $t_{r}\left(  n\right)  $ be the number of edges of
the $r$-partite Tur\'{a}n graph of order $n,$ we prove that:

(A) For all $r\geq2$ and all sufficiently small $\varepsilon>0,$ every graph
of sufficiently large order $n$ with $t_{r}\left(  n\right)  +1$ edges
contains a $K_{r}^{+}\left(  \left\lfloor c\ln n\right\rfloor ,\ldots
,\left\lfloor c\ln n\right\rfloor ,\left\lceil n^{1-\sqrt{c}}\right\rceil
\right)  .$

(B) For all $r\geq2,$ there exists $c>0$ such that every graph of sufficiently
large order $n$ with $t_{r}\left(  n\right)  +1$ edges contains a $K_{r}%
^{+}\left(  \left\lfloor c\ln n\right\rfloor ,\ldots,\left\lfloor c\ln
n\right\rfloor \right)  .$

These assertions extend results of Erd\H{o}s from 1963.

We also give corresponding stability results \medskip


\end{abstract}

\section{Introduction}

This note is part of an ongoing project aiming to renovate some classical
results in extremal graph theory, see, e.g., \cite{BoNi04},
\cite{Nik07,Nik07c}.

Let $t_{r}\left(  n\right)  $ be the number of edges of the $r$-partite
Tur\'{a}n graph of order $n.$ Tur\'{a}n's theorem implies that every graph on
$n$ vertices with $t_{r}\left(  n\right)  +1$ edges contains a $K_{r+1},$ the
complete graph of order $r+1.$ Thus, it is natural to ask:\medskip

\emph{Which supergraphs of }$K_{r+1}$\emph{ are present in graphs on }%
$n$\emph{ vertices with }$t_{r}\left(  n\right)  +1$\emph{ edges?}\medskip

A partial answer to this question was stated by Erd\H{o}s in \cite{Erd63} and
proved in \cite{ErSi73}, Theorem 1:\medskip

\emph{Let }$K_{r}^{+}\left(  s_{1},\ldots,s_{r}\right)  $\emph{ be the
complete }$r$\emph{-partite graph with parts of size }$s_{1}\geq2,s_{2}%
,\ldots,s_{r}$\emph{ with an edge added to the first part. For all }$s\geq
2,$\emph{ every graph of sufficiently large order }$n$\emph{ with }%
$t_{r}\left(  n\right)  +1$\emph{ edges contains a }$K_{r}^{+}\left(
s,\ldots,s\right)  .$\medskip

For $r=2,$ Erd\H{o}s \cite{Erd63} gave a stronger result:\medskip

\emph{For all sufficiently small }$\varepsilon>0,$\emph{ every graph of
sufficiently large order }$n$\emph{ with }$t_{2}\left(  n\right)  +1$\emph{
edges contains a }$K_{2}^{+}\left(  \left\lfloor c\ln n\right\rfloor
,\left\lceil n^{1-\varepsilon}\right\rceil \right)  $\emph{ for some }%
$c>0,$\emph{ independent of }$n.$\medskip

We extend these two results as follows.

\begin{theorem}
\label{Th1}Let $r\geq2,$ $2/\ln n\leq c\leq r^{-\left(  r+7\right)  \left(
r+1\right)  },$ and $G$ be a graph of order $n$. If $G$ has $t_{r}\left(
n\right)  +1$ edges, then $G$ contains a $K_{r}^{+}\left(  \left\lfloor c\ln
n\right\rfloor ,\ldots,\left\lfloor c\ln n\right\rfloor ,\left\lceil
n^{1-\sqrt{c}}\right\rceil \right)  .$
\end{theorem}

Here is a simpler version of this assertion.

\begin{theorem}
\label{Th1.1}Let $r\geq2,$ $c=r^{-\left(  r+7\right)  \left(  r+1\right)  },$
$n\geq e^{2/c},$ and $G$ be a graph of order $n$. If $G$ has $t_{r}\left(
n\right)  +1$ edges, then $G$ contains a $K_{r}^{+}\left(  \left\lfloor c\ln
n\right\rfloor ,\ldots,\left\lfloor c\ln n\right\rfloor \right)  .$
\end{theorem}

Theorems \ref{Th1} and \ref{Th1.1} have corresponding stability results.

\begin{theorem}
\label{Th2}Let $r\geq2,$ $2/\ln n\leq c\leq r^{-\left(  r+7\right)  \left(
r+1\right)  }/2,$ $0<\alpha<r^{-8}/8,$ and $G$ be a graph of order $n$. If $G$
has $\left\lceil \left(  1-1/r-\alpha\right)  n^{2}/2\right\rceil $ edges,
then $G$ satisfies one of the conditions:

(i) $G$ contains a $K_{r}^{+}\left(  \left\lfloor c\ln n\right\rfloor
,\ldots,\left\lfloor c\ln n\right\rfloor ,\left\lceil n^{1-2\sqrt{c}%
}\right\rceil \right)  ;$

(ii) $G$ contains an induced $r$-partite subgraph $G_{0}$ of order at least
$\left(  1-\sqrt{2\alpha}\right)  n$ and with minimum degree $\delta\left(
G_{0}\right)  >\left(  1-1/r-2\sqrt{2\alpha}\right)  n.$
\end{theorem}

\begin{theorem}
\label{Th2.1}Let $r\geq2,$ $c=r^{-\left(  r+7\right)  \left(  r+1\right)
}/2,$ $0<\alpha<r^{-8}/8,$ $n\geq e^{2/c},$ and $G$ be a graph of order $n$.
If $G$ has $\left\lceil \left(  1-1/r-\alpha\right)  n^{2}/2\right\rceil $
edges, then $G$ satisfies one of the conditions:

(i) $G$ contains a $K_{r}^{+}\left(  \left\lfloor c\ln n\right\rfloor
,\ldots,\left\lfloor c\ln n\right\rfloor \right)  ;$

(ii) $G$ contains an induced $r$-partite subgraph $G_{0}$ of order at least
$\left(  1-\sqrt{2\alpha}\right)  n$ and with minimum degree $\delta\left(
G_{0}\right)  >\left(  1-1/r-2\sqrt{2\alpha}\right)  n.$
\end{theorem}

\subsubsection*{Remarks}

\begin{itemize}
\item[-] The relations between $c$ and $n$ in Theorems \ref{Th1} and \ref{Th2}
need some explanation. First, for fixed $c,$ they show how large must be $n$
to get valid conclusions. But, in fact, the relations are subtler, for $c$
itself may depend on $n,$ e.g., letting $c=1/\ln\ln n,$ the conclusions are
meaningful for sufficiently large $n.$

\item[-] Note that, in Theorems \ref{Th1} and \ref{Th2}, if the conclusion
holds for some $c,$ it holds also for $0<c^{\prime}<c,$ provided $n$ is
sufficiently large.

\item[-] The stability conditions in Theorems \ref{Th2} and \ref{Th2.1} are
stronger than the conditions in the stability theorems of \cite{Erd68},
\cite{Sim68} and \cite{Nik07a}. Indeed, condition \emph{(ii)} implies that
$G_{0}$ is an induced, almost balanced, and almost complete $r$-partite graph
containing almost all the vertices of $G;$

\item[-] The exponents $1-\sqrt{c}$ and $1-2\sqrt{c}$ in Theorems \ref{Th1}
and \ref{Th2} are far from the best ones, but are simple.
\end{itemize}

In our proofs we apply the following two technical statements, which may be
useful elsewhere.

\begin{lemma}
\label{le1} Let $0<\alpha\leq1,$ $1\leq c\ln n\leq\alpha m/2+1,$ and let $F$
be a bipartite graph with parts $A$ and $B$ of size $m$ and $n.$ If $e\left(
F\right)  \geq\alpha mn,$ then $F$ contains a $K_{2}\left(  s,t\right)  $ with
parts $S\subset A$ and $T\subset B$ such that $\left\vert S\right\vert
=\left\lfloor c\ln n\right\rfloor $ and $\left\vert T\right\vert
=t>n^{1-c\ln\alpha/2}$.
\end{lemma}

\begin{theorem}
\label{thv3}Let $r\geq2,$ $2/\ln n\leq c\leq r^{-\left(  r+8\right)  r}$ and
$G$ be a graph $G$ of order $n.$ If $G$ contains a $K_{r+1}$ and has minimum
degree $\delta\left(  G\right)  >\left(  1-1/r-1/r^{4}\right)  n,$ then $G$
contains a
\[
K_{r}^{+}\left(  \left\lfloor c\ln n\right\rfloor ,\ldots,\left\lfloor c\ln
n\right\rfloor ,\left\lceil n^{1-cr^{3}}\right\rceil \right)  .
\]

\end{theorem}

The next section contains notation and results needed to prove the theorems.
The proofs are presented in Section \ref{pf}.

\section{Preliminary results}

Our notation follows \cite{Bol98}; given a graph $G,$ we write:$\smallskip$

- $V\left(  G\right)  $ for the vertex set of $G$ and $\left\vert G\right\vert
$ for $\left\vert V\left(  G\right)  \right\vert ;$

- $E\left(  G\right)  $ for the edge set of $G$ and $e\left(  G\right)  $ for
$\left\vert E\left(  G\right)  \right\vert ;$

- $\Gamma\left(  u\right)  $ for the set of neighbors of a vertex $u$ and
$d\left(  u\right)  $ for $\left\vert \Gamma\left(  u\right)  \right\vert ;$

- $\delta\left(  G\right)  $ for the minimum degree of $G;$

- $G\left[  U\right]  $ for the subgraph of $G$ induced by a set $U\subset
V\left(  G\right)  ;$

- $H+u$ for $G\left[  V\left(  H\right)  \cup\left\{  u\right\}  \right]  ,$
where $H\subset G$ is a subgraph and $u\in V\left(  G\right)  ;$

- $K_{r}\left(  G\right)  $ for the set of $r$-cliques of $G$ and
$k_{r}\left(  G\right)  $ for $\left\vert K_{r}\left(  G\right)  \right\vert
;$

- $K_{s}\left(  M\right)  $ for the set of $s$-cliques contained in members of
a set $M\subset K_{r}\left(  G\right)  ;$

- $K_{r}\left(  s_{1},\ldots,s_{r}\right)  $ for the complete $r$-partite
graph with parts of size $s_{1},\ldots,s_{r}.$\bigskip

An\emph{ }$r$\emph{-joint }of size $t$ is the union of $t$ distinct
$r$-cliques sharing an edge. Write $js_{r}\left(  G\right)  $ for the maximum
size of an $r$-joint in $G.\medskip$

Given a set $M\subset K_{r}\left(  G\right)  ,$ we say that $M$ \emph{covers}
a subgraph $H\subset G,$ if $E\left(  H\right)  \subset K_{2}\left(  M\right)
$.\bigskip

The following facts play crucial roles in our proofs.

\begin{fact}
[\cite{BoNi04}, Lemma 1]\label{leNSMM}Let $r\geq2$ and $c\geq0,$ and $G$ be a
graph of order $n.$ If%
\[
e\left(  G\right)  >\left(  1-1/r+c\right)  n^{2}/2,
\]
then%
\[
k_{r+1}\left(  G\right)  >c\frac{r^{2}}{r+1}\left(  \frac{n}{r}\right)
^{r+1}.
\]
$\hfill\square$
\end{fact}

\begin{fact}
[\cite{BoNi04}, Lemma 6]\label{leKd} Let $r\geq2,$ and $G$ be a graph of order
$n.$ If $G$ contains a $K_{r+1}$ and $\delta\left(  G\right)  >\left(
1-1/r-1/r^{4}\right)  n,$ then $js_{r+1}\left(  G\right)  >n^{r-1}%
/r^{r+3}.\hfill\square$
\end{fact}

\begin{fact}
[\cite{BoNi04}, Theorem 7]\label{Thexj}Let $r\geq2,$ $n>r^{8}$, and $G$ be a
graph of order $n.$ If $e\left(  G\right)  >t_{r}\left(  n\right)  ,$ then $G$
has an induced subgraph $G^{\prime}$ of order $n^{\prime}>\left(
1-1/r^{2}\right)  n$\ such that either
\begin{equation}
e\left(  G^{\prime}\right)  >\left(  \frac{r-1}{2r}+\frac{1}{r^{4}\left(
r^{2}-1\right)  }\right)  \left(  n^{\prime}\right)  ^{2} \label{prop4}%
\end{equation}
or%
\begin{equation}
K_{r+1}\subset G^{\prime}\text{, \ \ \ and \ \ \ }\delta\left(  G^{\prime
}\right)  >\left(  1-1/r-1/r^{4}\right)  n^{\prime}. \label{prop3}%
\end{equation}
$\hfill\square$
\end{fact}

\begin{fact}
[\cite{Nik07}, Theorem 1]\label{ES}Let $r\geq2,$ $\alpha^{r}\ln n\geq1,$ and
$G$ be a graph of order $n$. Every set $M\subset K_{r}\left(  G\right)  $
satisfying $\left\vert M\right\vert \geq\alpha n^{r}$ covers a $K_{r}\left(
s,\ldots s,t\right)  $ with $s=\left\lfloor \alpha^{r}\ln n\right\rfloor $ and
$t>n^{1-\alpha^{r-1}}.\hfill\square$
\end{fact}

\section{\label{pf}Proofs}

\begin{proof}
[\textbf{Proof of Lemma \ref{le1}}]Set $s=\left\lfloor c\ln n\right\rfloor $
and let
\[
t=\max\left\{  x:\text{there exists }K_{2}\left(  s,x\right)  \subset F\text{
with part of size }s\text{ in }A\right\}  .
\]
Thus $d\left(  X\right)  \leq t$ for each $X\subset A$ with $\left\vert
X\right\vert =s,$ and so,%
\begin{equation}
t\binom{m}{s}\geq%
{\textstyle\sum\limits_{X\subset A,\left\vert X\right\vert =s}}
d\left(  X\right)  =%
{\textstyle\sum\limits_{u\in B}}
\binom{d\left(  u\right)  }{s}. \label{in1}%
\end{equation}
Setting
\[
f\left(  x\right)  =\left\{
\begin{array}
[c]{cc}%
\binom{x}{s} & \text{if }x\geq s-1\\
0 & \text{if }x<s-1,
\end{array}
\right.
\]
and noting that $f\left(  x\right)  $ is a convex function, we find that,%
\[%
{\textstyle\sum\limits_{u\in B}}
\binom{d\left(  u\right)  }{s}=%
{\textstyle\sum\limits_{u\in B}}
f\left(  d\left(  u\right)  \right)  \geq nf\left(  \frac{1}{n}%
{\textstyle\sum\limits_{u\in B}}
d\left(  u\right)  \right)  =n\binom{e\left(  F\right)  /n}{s}\geq
n\binom{\alpha m}{s}.
\]
Combining this inequality with (\ref{in1}) and rearranging, we find that%
\[
t\geq n\frac{\alpha m\left(  cm-1\right)  \cdots\left(  \alpha m-s+1\right)
}{m\left(  m-1\right)  \cdots\left(  m-s+1\right)  }>n\left(  \frac{\alpha
m-s+1}{m}\right)  ^{s}\geq n\left(  \frac{\alpha}{2}\right)  ^{s}\geq
n^{1+c\ln\left(  \alpha/2\right)  },
\]
completing the proof.
\end{proof}

\begin{proof}
[\textbf{Proof of Theorem \ref{thv3}}]Let $r,c,n,$ and the graph $G$ satisfy
the conditions of the theorem. Note first that for every $R\in K_{r-1}\left(
G\right)  ,$
\begin{align}
d\left(  R\right)   &  =\left\vert
{\textstyle\bigcap\limits_{u\in R}}
\Gamma\left(  u\right)  \right\vert \geq%
{\textstyle\sum\limits_{u\in R}}
d\left(  u\right)  -\left(  r-2\right)  n\geq\left(  r-1\right)  \delta\left(
G\right)  -\left(  r-2\right)  n\nonumber\\
&  >\left(  \left(  r-1\right)  \left(  \left(  1-\frac{1}{r}\right)
-\frac{1}{r^{4}}\right)  -\left(  r-2\right)  \right)  n>\frac{n}{r^{2}}.
\label{c1}%
\end{align}

Also, Fact \ref{leKd} implies that%
\[
js_{r+1}\left(  G\right)  >\frac{n^{r-1}}{r^{r+3}}>\left(  1-\frac{1}{r^{2}%
}\right)  ^{r-1}\frac{n^{r-1}}{r^{r+3}}>\left(  1-\frac{r-1}{r^{2}}\right)
\frac{n^{r-1}}{r^{r+3}}>\frac{n^{r-1}}{r^{r+4}}.
\]
Thus, there exists an edge $uv\in E\left(  G\right)  $ contained in more than
$n^{r-1}/r^{r+4}$ distinct $\left(  r+1\right)  $-cliques of $G.$ Letting
$B=\Gamma\left(  u\right)  \cap\Gamma\left(  v\right)  \cap V\left(  G\right)
,$ we see that
\begin{equation}
k_{r-1}\left(  G\left[  B\right]  \right)  >n^{r-1}/r^{r+4}. \label{c2}%
\end{equation}

Define the set $X$ as%
\[
X=\left\{  R:R\in K_{r}\left(  G\right)  \text{ \ and \ }\left\vert R\cap
B\right\vert \geq r-1\right\}  .
\]
In view of (\ref{c1}) and (\ref{c2}), we find that
\[
\left\vert X\right\vert \geq\frac{1}{r}%
{\displaystyle\sum_{P\in K_{r-1}\left(  G\left[  B\right]  \right)  }}
d\left(  P\right)  >\frac{1}{r}\times\frac{n}{r^{2}}\times\frac{n^{r-1}%
}{r^{r+4}}=\frac{n^{r}}{r^{r+7}},
\]

For a set $N\subset K_{r}\left(  G\right)  $ and a clique $R\in K_{r-1}\left(
N\right)  $ let $d_{N}\left(  R\right)  $ be the number of members of $N$
containing $R$. We claim that there exists $Y\subset X$ with $\left\vert
Y\right\vert >n^{r}/r^{r+8}$ such that $d_{Y}\left(  R\right)  >n/r^{r+8}$ for
all $R\in K_{r-1}\left(  Y\right)  .$ Indeed, set $Y=X$ and apply the
following procedure:\medskip

\textbf{While}\emph{ there exists an }$R\in K_{r-1}\left(  Y\right)  $\emph{
with }$d_{Y}\left(  R\right)  \leq n/r^{r+8}$ \textbf{do}\emph{ }

\qquad\emph{Remove from }$Y$\emph{ all }$r$\emph{-cliques containing }%
$R$\emph{.}\medskip

When the procedure stops, $d_{Y}\left(  R\right)  >n/r^{r+8}$ for all $R\in
K_{r-1}\left(  Y\right)  ,$ and
\[
\left\vert X\right\vert -\left\vert Y\right\vert \leq\left\vert K_{r-1}\left(
X\right)  \right\vert \frac{n}{r^{r+8}}\leq\binom{n}{r-1}\frac{n}{r^{r+8}%
}<\frac{1}{r^{r+8}}n^{r},
\]
implying that $\left\vert Y\right\vert >n^{r}/r^{r+8},$ as claimed.

Since
\[
\left\vert K_{r-1}\left(  Y\right)  \right\vert \geq r\left\vert Y\right\vert
/n>r\times r^{-r-8}n^{r}/n=n^{r-1}/r^{r+7},
\]
by Fact \ref{ES}, $K_{r-1}\left(  Y\right)  $ covers a subgraph $H=K_{r-1}%
\left(  m,\ldots,m\right)  $ with $m=\left\lfloor r^{-\left(  r+7\right)
\left(  r-1\right)  }\ln n\right\rfloor .$

Select a set $A$ of $m$ disjoint $\left(  r-1\right)  $-cliques in $H$ and
define a bipartite graph $F$ with parts $A$ and $B,$ joining $R\in A$ to $v\in
B$ if $R+v\in Y.$

Let $\alpha=1/r^{r+8}$ and set $s=\left\lfloor c\ln n\right\rfloor .$ Since
\[
d_{Y}\left(  R\right)  >\frac{1}{r^{r+8}}n\geq\alpha n
\]
for all $R\in K_{r-1}\left(  Y\right)  ,$ we have $e\left(  F\right)  >\alpha
mn.$ Also, we find that%
\[
s\leq c\ln n\leq\frac{1}{r^{\left(  r+8\right)  r}}\ln n\leq\frac{1}{2r^{r+8}%
}\times\frac{1}{r^{\left(  r+7\right)  \left(  r-1\right)  }}\ln n\leq
\frac{\alpha}{2}m+1.
\]

Hence, by Fact \ref{le1}, $H$ contains a $K_{2}\left(  s,t\right)  $ with
parts $S\subset A$ and $T\subset B$ such that $\left\vert S\right\vert =s$ and
$\left\vert T\right\vert =t>n^{1-c\ln\alpha/2}$. A routine calculation shows
that for $r\geq2,$
\[
\ln\alpha/2=\ln\frac{1}{2r^{r+8}}\geq-r^{3},
\]
and so, $t>n^{1-cr^{3}}.$

Letting $H^{\ast}$ be the subgraph of $H$ induced by the union of the members
of $S,$ we see that $H^{\ast}=K_{r-1}\left(  s,\ldots,s\right)  $. Note that
at least $\left(  r-2\right)  $ of the parts of $H^{\ast}$ belong to $B,$ for
otherwise we can select an $\left(  r-1\right)  $-clique $Q$ in $H^{\ast}$
with two vertices outside $B,$ and so, every $R\in Y$ containing $Q$ has two
vertices outside $B.$ This is a contradiction since $Y\subset X$ and all
members of $X$ intersect $B$ in at least $r-1$ vertices.

Let $H_{1},\ldots,H_{r-1}$ be the parts of $H^{\ast},$ and assume by symmetry
that $H_{i}\subset B$ for $i=2,\ldots,r-1.$ Remove two vertices from $H_{1},$
add $u$ and $v$ to $H_{1},$ and write $H_{1}^{\prime}$ for the resulting set.
Clearly the sets $H_{1}^{\prime},H_{2}\ldots,H_{r-1},T$ induce a subgraph
containing a $K_{r}^{+}\left(  \left\lfloor c\ln n\right\rfloor ,\ldots
,\left\lfloor c\ln n\right\rfloor ,\left\lceil n^{1-cr^{3}}\right\rceil
\right)  ,$ completing the proof.\bigskip
\end{proof}

\begin{proof}
[\textbf{Proof of Theorem \ref{Th1}}]Let $G$ be a graph of order $n$ with
$t_{r}\left(  n\right)  +1$ edges. Fact \ref{Thexj} implies that there exists
an induced subgraph $G^{\prime}\subset G$ of order $n^{\prime}>\left(
1-1/r^{2}\right)  n$\ such that either (\ref{prop4}) or (\ref{prop3}) holds.

Assume first that $G^{\prime}$ satisfies condition (\ref{prop4}). Fact
\ref{leNSMM} implies that%
\begin{align*}
k_{r+1}\left(  G\right)   &  \geq k_{r+1}\left(  G^{\prime}\right)  >\frac
{2}{r^{4}\left(  r^{2}-1\right)  }\times\frac{r^{2}}{r+1}\times\left(
\frac{n^{\prime}}{r}\right)  ^{r+1}\\
&  >\frac{2}{r^{2}\left(  r^{2}-1\right)  \left(  r+1\right)  }\times\left(
1-\frac{1}{r^{2}}\right)  ^{r+1}\times\left(  \frac{n}{r}\right)  ^{r+1}\\
&  >\frac{2}{r^{2}\left(  r^{2}-1\right)  \left(  r+1\right)  }\times\left(
1-\frac{r+1}{r^{2}}\right)  \times\left(  \frac{n}{r}\right)  ^{r+1}\\
&  >\frac{2\left(  r^{2}-r-1\right)  }{r^{4}\left(  r^{2}-1\right)  \left(
r+1\right)  }\times\left(  \frac{n}{r}\right)  ^{r+1}>\frac{1}{r^{r+7}}%
n^{r+1}>c^{1/\left(  r+1\right)  }n^{r+1}.
\end{align*}
Hence, by Fact \ref{ES}, $G$ contains a $K_{r+1}\left(  s,\ldots,s,t\right)  $
with $s=\left\lfloor c\ln n\right\rfloor $ and
\[
t>n^{1-c^{r/\left(  r+1\right)  }}>n^{1-\sqrt{c}}.
\]
Then, obviously, $G$ contains a $K_{r}^{+}\left(  \left\lfloor c\ln
n\right\rfloor ,\ldots,\left\lfloor c\ln n\right\rfloor ,\left\lceil
n^{1-\sqrt{c}}\right\rceil \right)  ,$ completing the proof.\bigskip

Finally, assume that $G^{\prime}$ satisfies condition (\ref{prop3}). Applying
Theorem \ref{thv3}, we see that $G^{\prime}$ contains a%
\[
K_{r}^{+}\left(  \left\lfloor 2c\ln n^{\prime}\right\rfloor ,\ldots
,\left\lfloor 2c\ln n^{\prime}\right\rfloor ,\left\lceil \left(  n^{\prime
}\right)  ^{1-2cr^{3}}\right\rceil \right)  .
\]
To complete the proof, note that
\[
2c\ln n^{\prime}\geq2c\ln\left(  1-\frac{1}{r^{2}}\right)  n\geq2\ln\left(
1-\frac{1}{r^{2}}\right)  +2\ln n\geq c\ln n
\]
and
\[
\left(  n^{\prime}\right)  ^{1-2cr^{3}}\geq\left(  1-\frac{1}{r^{2}}\right)
^{1-2cr^{3}}n^{1-2cr^{3}}\geq\left(  1-\frac{1}{r^{2}}\right)  n^{1-2cr^{3}%
}>n^{1-\sqrt{c}}.
\]

\end{proof}

\bigskip

\begin{proof}
[\textbf{Proof of Theorem \ref{Th2}}]Let $G$ be a graph of order $n$ with
$e\left(  G\right)  >\left(  1-1/r-\alpha\right)  n^{2}/2.$ Set $V=V\left(
G\right)  ,$ $\varepsilon=\sqrt{2\alpha},$ and define the set $M_{\varepsilon
}$ as%
\[
M_{\varepsilon}=\left\{  u\in V\left(  G\right)  :d\left(  u\right)
\leq\left(  1-1/r-\varepsilon\right)  n\right\}  .
\]

Assume that condition \emph{(i)} fails. We shall show that: \emph{(a)
}$\left\vert M_{\varepsilon}\right\vert <\varepsilon n;$ \emph{(b)} the graph
$G_{0}=G\left[  V\backslash M_{\varepsilon}\right]  $ satisfies condition
\emph{(ii)}.\bigskip

\emph{(a) The set }$M_{\varepsilon}$\emph{ satisfies }$\left\vert
M_{\varepsilon}\right\vert <\varepsilon n$\bigskip

Assume for a contradiction that $\left\vert M_{\varepsilon}\right\vert
\geq\varepsilon n$, select $M^{\prime}\subset M_{\varepsilon}$ with
\begin{equation}
\left\vert M^{\prime}\right\vert =\left\lfloor \varepsilon n\right\rfloor
\label{bnds}%
\end{equation}
and note that $M^{\prime}$ is nonempty since $\varepsilon n=\sqrt{2\alpha
}n>1.$ Letting $G^{\prime}=G\left[  V\backslash M^{\prime}\right]  ,$ we see
that%
\begin{align*}
e\left(  G\right)   &  =e\left(  G^{\prime}\right)  +e\left(  M^{\prime
},V\backslash M^{\prime}\right)  +e\left(  M^{\prime}\right)  \leq e\left(
G^{\prime}\right)  +%
{\textstyle\sum\limits_{u\in M^{\prime}}}
d\left(  u\right) \\
&  \leq e\left(  G^{\prime}\right)  +\left\vert M^{\prime}\right\vert \left(
1-1/r-\varepsilon\right)  n.
\end{align*}

Assume for a contradiction that
\[
e\left(  G^{\prime}\right)  >\frac{r-1}{2r}\left(  n-\left\vert M^{\prime
}\right\vert \right)  ^{2}%
\]
and set $p=n-\left\vert M^{\prime}\right\vert .$ In view of (\ref{bnds}), we
have
\[
p\geq n-\varepsilon n=\left(  1-\sqrt{2\alpha}\right)  n.
\]
Hence, by Theorem \ref{Th1}, $G$ contains a $K_{r}^{+}\left(  \left\lfloor
2c\ln p\right\rfloor ,\ldots,\left\lfloor 2c\ln p\right\rfloor ,\left\lceil
p^{1-\sqrt{2c}}\right\rceil \right)  .$ Since%
\[
2c\ln p\geq2c\ln\left(  1-\sqrt{2\alpha}\right)  n\geq2c\ln\left(  1-\frac
{1}{4r^{4}}\right)  n\geq c\ln n
\]
and%
\[
p^{1-\sqrt{2c}}\geq\left(  1-\sqrt{2\alpha}\right)  ^{1-\sqrt{2c}}%
n^{1-\sqrt{2c}}>\left(  1-\sqrt{2\alpha}\right)  n^{1-\sqrt{2c}}%
>n^{1-2\sqrt{c}},
\]
this contradicts the assumption that \emph{(i)} fails.

Hereafter, we assume that%
\[
e\left(  G^{\prime}\right)  \leq\frac{r-1}{2r}\left(  n-\left\vert M^{\prime
}\right\vert \right)  ^{2}.
\]
From
\[
e\left(  G^{\prime}\right)  \geq e\left(  G\right)  -%
{\textstyle\sum\limits_{u\in M}}
d\left(  u\right)  \geq\left(  1-1/r-\alpha\right)  n^{2}/2-\left\vert
M^{\prime}\right\vert \left(  1-1/r-\varepsilon\right)  n,
\]
we obtain
\[
\frac{r-1}{2r}\left(  n-\left\vert M^{\prime}\right\vert \right)  ^{2}%
\geq\left(  \frac{r-1}{r}-\alpha\right)  \frac{n^{2}}{2}-\left\vert M^{\prime
}\right\vert \left(  \frac{r-1}{r}-\varepsilon\right)  n.
\]
After some algebra, we find that%
\[
\left\vert M^{\prime}\right\vert <\left(  \varepsilon-\sqrt{\varepsilon
^{2}-\alpha}\right)  n=\varepsilon\left(  1-\sqrt{1/2}\right)  n
\]
or%
\[
\left\vert M^{\prime}\right\vert >\left(  \varepsilon+\sqrt{\varepsilon
^{2}-\alpha}\right)  n=\varepsilon\left(  1+\sqrt{1/2}\right)  n,
\]
contradicting (\ref{bnds}) in view of $\varepsilon\sqrt{1/2}n=\sqrt{2\alpha
}n>\sqrt{2}.$ Therefore, $\left\vert M_{\varepsilon}\right\vert <\varepsilon
n.$\bigskip

\emph{(b) The graph }$G_{0}=G\left[  V\backslash M_{\varepsilon}\right]
$\emph{ satisfies condition (ii).}\bigskip

By our choice of $M_{\varepsilon},$ for every $u\in V\backslash M_{\varepsilon
},$ we have $d\left(  u\right)  >\left(  1-1/r-\varepsilon\right)  n;$ thus
\[
\delta\left(  G_{0}\right)  >\left(  1-1/r-\varepsilon\right)  n-\left\vert
M_{\varepsilon}\right\vert >\left(  1-1/r-2\varepsilon\right)  n=\left(
1-1/r-2\sqrt{2\alpha}\right)  n,
\]
and so, $\delta\left(  G_{0}\right)  $ satisfies the required condition. All
that remains to prove is that $G_{0}$ is $r$-partite.

If $G_{0}$ contains a $K_{r+1}$, in view of
\[
\delta\left(  G_{0}\right)  >\left(  1-1/r-2\sqrt{2\alpha}\right)  n>\left(
1-1/r-1/r^{4}\right)  n,
\]
using Theorem \ref{thv3} as in the proof of Theorem \ref{Th1}, we see that $G$
contains a
\[
K_{r}^{+}\left(  \left\lfloor c\ln n\right\rfloor ,\ldots,\left\lfloor c\ln
n\right\rfloor ,\left\lceil n^{1-\sqrt{c}}\right\rceil \right)  ,
\]
contradicting our assumption. Thus, $G_{0}$ is $K_{r+1}$-free. In view of
\[
\delta\left(  G_{0}\right)  >\left(  1-1/r-1/r^{4}\right)  n>\left(
1-\frac{3}{3r-1}\right)  \left\vert G_{0}\right\vert ,
\]
the theorem of Andr\'{a}sfai, Erd\H{o}s and S\'{o}s \cite{AES74} implies that
$G_{0}$ is $r$-partite, completing the proof.
\end{proof}

We omit the proofs of Theorems \ref{Th1.1} and \ref{Th2.1}, since they are
easy consequences of Theorem \ref{Th1} and \ref{Th2}.

\subsubsection*{Concluding remark}

Finally, a word about the project mentioned in the introduction: in this
project we aim to give wide-range results that can be used further, adding
more integrity to extremal graph theory.

\end{document}